\documentclass{article}
\usepackage{amssymb,amsmath,graphicx,ucs,hyperref,fullpage}
\usepackage[utf8x]{inputenc}

\title{Compact Packings are not always the Densest}

\author{
Thomas Fernique\footnote{LIPN, Univ. Paris Nord \& CNRS, 93430 Villetaneuse, France.}
\and Daria Pchelina\footnote{LIPN, Univ. Paris Nord \& CNRS, 93430 Villetaneuse, France.}
}
\date{}

\begin{document}
\maketitle

A {\em disc packing} is a union of disjoint interior discs in the Euclidean plane.
It is said to be {\em compact} if its {\em contact graph}, i.e., the graph which connects the center of adjacent discs, is triangulated.

There is only one compact packing by unit discs, called the {\em hexagonal compact packing}.
It is proven in \cite{FT43} to have maximal density among all the packings by unit discs.

There are exactly $9$ values $r<1$ which allow compact packings by discs of sizes $1$ and $r$ (where the two disc sizes appear) \cite{Ken06}.
In each case the maximal density has been proven to be reached for a periodic compact packing \cite{BF20} (see also \cite{Hep00,Hep03,Ken04}).

\begin{figure}[hbt]
\centering
\begin{tabular}{lll}
  \includegraphics[width=0.3\textwidth]{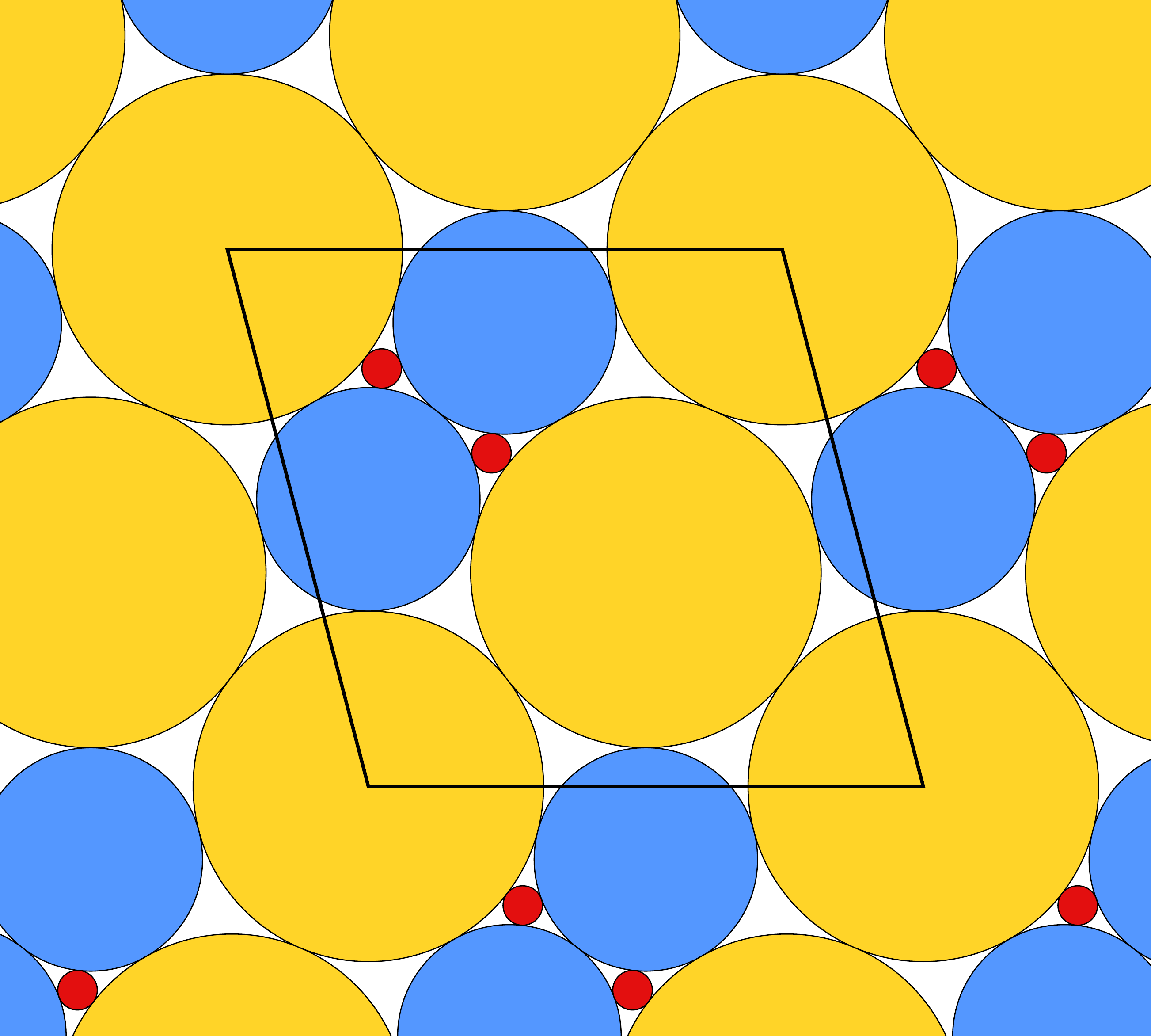} &
  \includegraphics[width=0.3\textwidth]{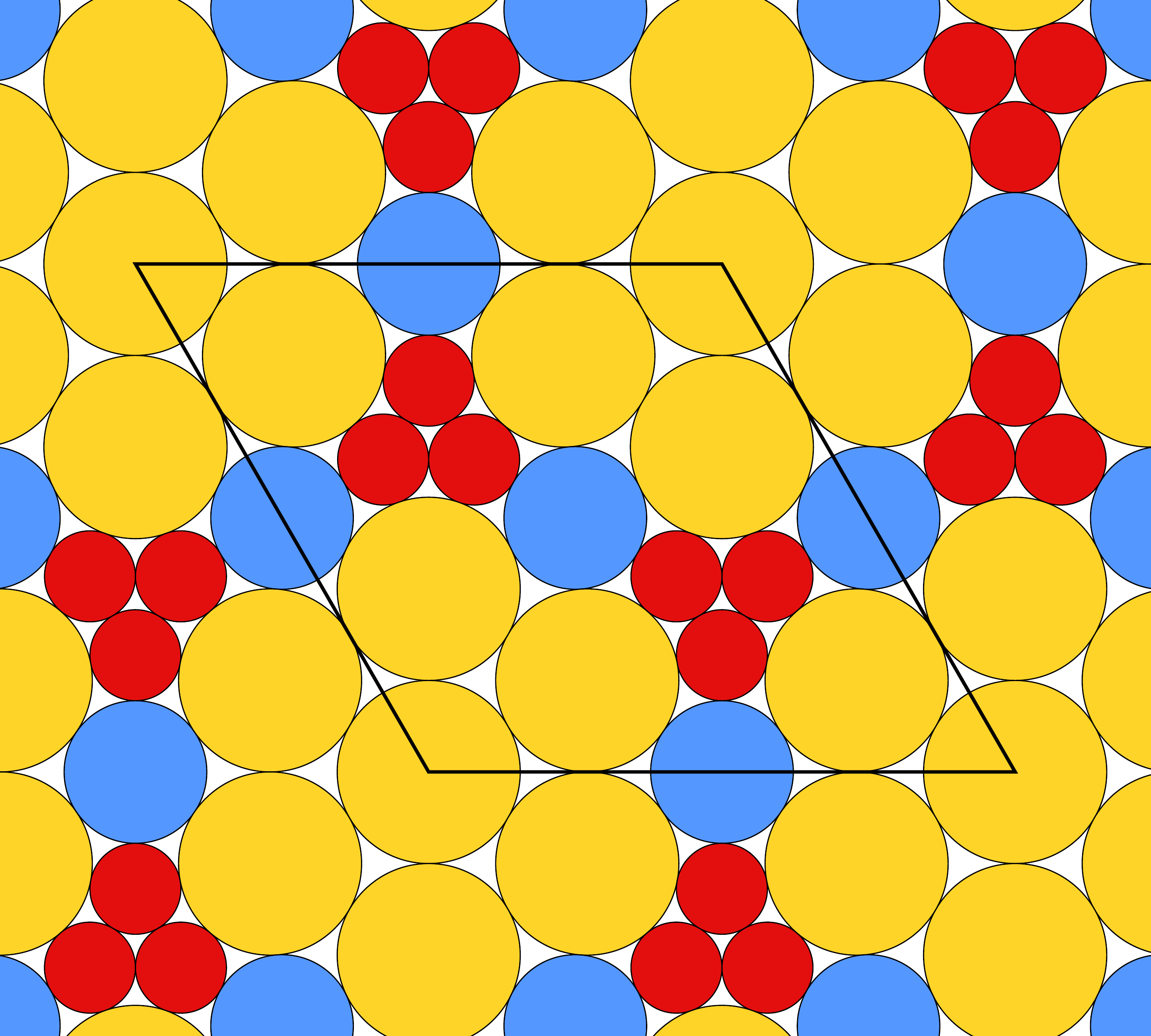}&
  \includegraphics[width=0.3\textwidth]{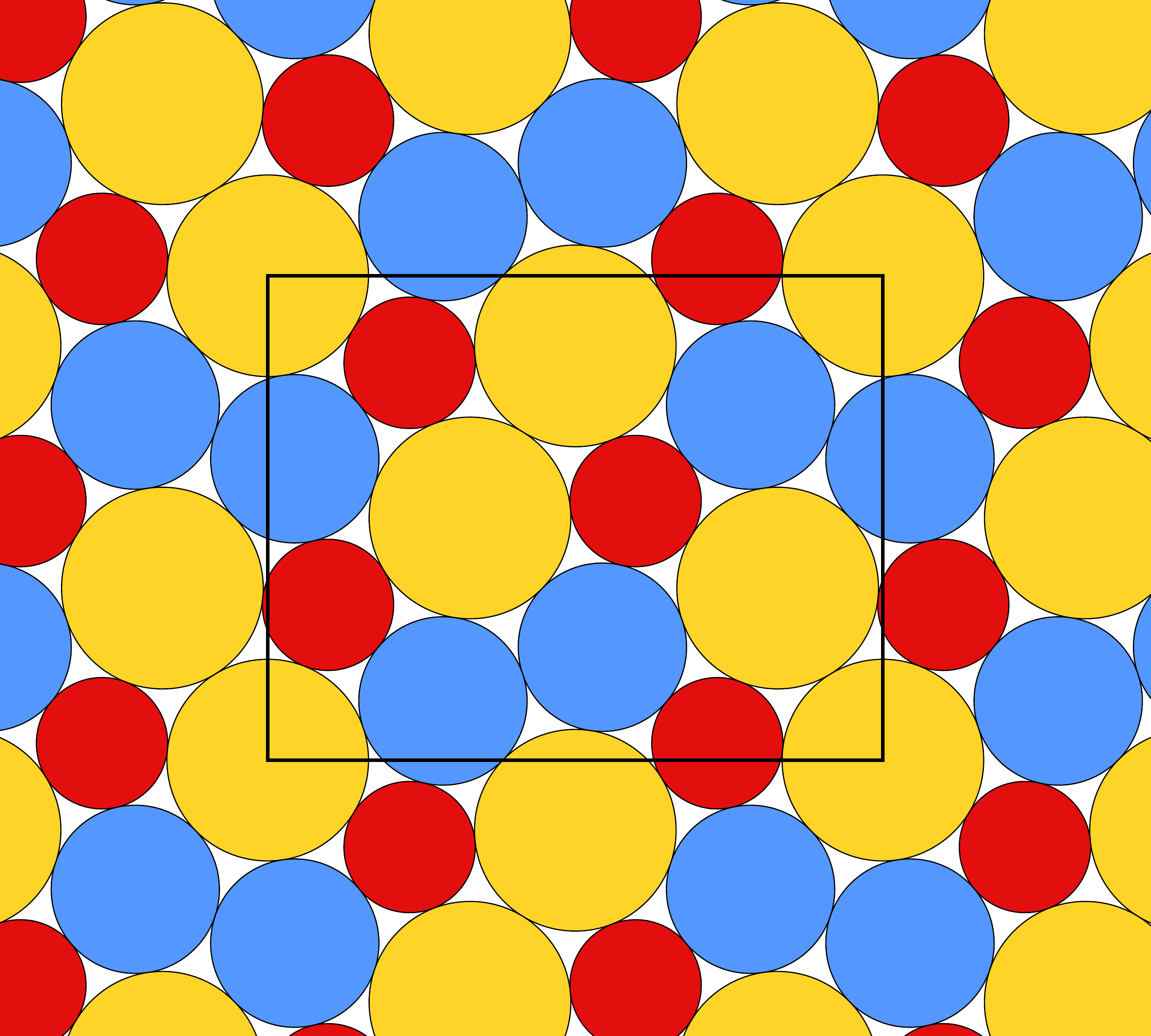} \\
\end{tabular}
\caption{
Three ternary compact packings numbered $28$, $110$ and $53$ in \cite{FHS20}.
}
\label{fig:conj_connelly}
\end{figure}

There are exactly $164$ values $s<r<1$ which allow compact packings by discs of sizes $1$, $r$ and $s$ (where the three disc sizes appear) \cite{FHS20}.
However, the maximal density is {\em not} always reached for a compact packing.
A simple counterexample is depicted on Fig.~\ref{fig:conj_connelly}, left.
In this case, the smallest discs can be inserted in each holes between two large and one medium discs, yielding a noncompact packing which is more dense.
This counterexample can be ruled out by considering only {\em saturated} packings, that is, such that no further disc can be added.
In \cite{CGSY18}, it has been conjectured that if the densest compact packing is saturated, then the maximal density is reached for a compact packing.

\begin{figure}[hbt]
\centering
\includegraphics[width=0.25\textwidth]{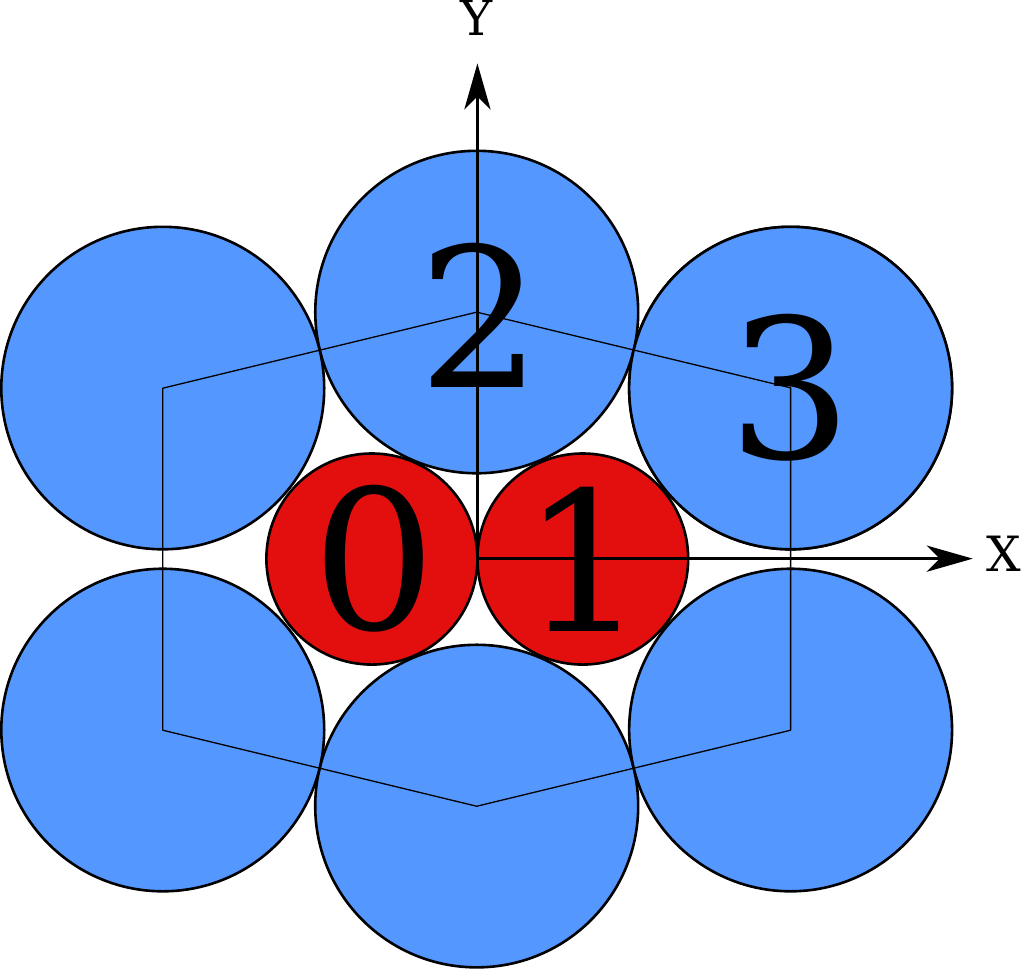}
\caption{Fundamental domain of a noncompact packing by the small and medium discs of the compact packing depicted in Fig.~\ref{fig:conj_connelly}, middle.
}
\label{fig:110}
\end{figure}

Unfortunately, the packing depicted in Fig~\ref{fig:conj_connelly}, middle, yields a new counterexample: it is saturated but the discs can form a noncompact packing which is more dense.
Let us give details.
The discs have size $1$, $r\approx 0.779$ and $s\approx 0.497$, where $r$ and $s$ are roots of, respectively:
$$
9x^8 + 12x^7 - 242x^6 + 436x^5 + 665x^4 - 2680x^3 + 2680x^2- 1056x + 144,
$$
$$
9x^8 - 120x^7 - 380x^6 + 2056x^5 + 12846x^4 - 29672x^3 + 15220x^2 - 2088x + 81.
$$
The point is that the ratio $q:=s/r\approx 0.6378$ is very close to one of the $9$ ratios which allow a compact packing by two sizes of discs (namely $r_1\approx 0.6375$).
It actually allows an ``almost compact'' packing, see Fig.~\ref{fig:110}.

In Fig.~\ref{fig:110}, a cartesian coordinate system is depicted and discs have been numbered.
Disc $0$ is centered in $(-q,0)$, disc $1$ in $(q,0)$ and disc 2 in $(0,\sqrt{2q+1})$.
A computation shows that disc $3$ has Y-coordinate:
$$
\frac{2 q\sqrt{q(q + 2)} + {\left(q^{2} + 2 q - 1\right)} \sqrt{2 q + 1}}{q^{2} + 2 q + 1}>1.0007.
$$
Further computations (see joined SageMath code \verb+110.sage+ for full details) shows that the density of this packing is greater than $0.9105$, while the density of the packing depicted in Fig.~\ref{fig:110}, middle, is less than $0.9104$.
This yields the claimed counterexample.
It uses only two of the three disc sizes, but it is easy to remove a positive but small enough proportion of discs so that the density is still larger than the one of the packing depicted in Fig.~\ref{fig:110}, middle, but the removed discs create holes large enough to contain each a large disc of radius $1$.

To conclude, let us mention the rightmost packing in Fig.~\ref{fig:conj_connelly}.
This one has indeed been proven in \cite{Fer19b} to maximize the density.
The maximal density of many of the $164$ cases is still open.
For which ones the conjecture stated in \cite{CGSY18} holds?

\bibliographystyle{alpha}
\bibliography{conjecture_connelly.bib}
\end{document}